\magnification=1200

\def\Q{{\bf {Q}}}
\def\K{{\bf {K}}}
  
\def\N{{\bf N}} 
     
\def\R{{\bf R}}

\def\C{{\bf C}}

\def\og{\leavevmode\raise.3ex\hbox{$\scriptscriptstyle 
\langle\!\langle\,$}}
\def \fg {\leavevmode\raise.3ex\hbox{$\scriptscriptstyle 
\rangle\!\rangle\,\,$}}

\def\house#1{\setbox1=\hbox{$\,#1\,$}%
\dimen1=\ht1 \advance\dimen1 by 2pt \dimen2=\dp1 \advance\dimen2 by 2pt
\setbox1=\hbox{\vrule height\dimen1 depth\dimen2\box1\vrule}%
\setbox1=\vbox{\hrule\box1}%
\advance\dimen1 by .4pt \ht1=\dimen1
\advance\dimen2 by .4pt \dp1=\dimen2 \box1\relax}

  \def\eps{{\varepsilon}}

\def\sm{\smallskip}  \def\noi{\noindent}

\def\build#1_#2^#3{\mathrel{\mathop{\kern 0pt#1}\limits_{#2}^{#3}}}

\def\date {le\ {\the\day}\ \ifcase\month\or janvier
\or fevrier\or mars\or avril\or mai\or juin\or juillet\or
ao\^ut\or septembre\or octobre\or novembre
\or d\'ecembre\fi\ {\oldstyle\the\year}}

\font\fivegoth=eufm5 \font\sevengoth=eufm7 \font\tengoth=eufm10

\newfam\gothfam \scriptscriptfont\gothfam=\fivegoth
\textfont\gothfam=\tengoth \scriptfont\gothfam=\sevengoth
\def\goth{\fam\gothfam\tengoth}

\def\pp{{\goth p}}

\def\smallsquare{\vbox{\hrule\hbox{\vrule height 1 ex\kern 1 ex\vrule}\hrule}}
\def\cqfd{\hfill \smallsquare\vskip 3mm}


\centerline{}

\vskip 4mm

\centerline{
\bf On the complexity of algebraic numbers I.
Expansions in integer bases}

\vskip 8mm
\centerline{Boris A{\sevenrm DAMCZEWSKI} \ (Lyon) \
\& \ Yann B{\sevenrm UGEAUD} \footnote{}{\rm 
} (Strasbourg)}

{\narrower\narrower
\vskip 12mm

\proclaim Abstract. {
Let $b \ge 2$ be an integer. We prove that the $b$-adic
expansion of every irrational algebraic number cannot have low complexity.
Furthermore, we establish that irrational morphic numbers are
transcendental, for a wide class of morphisms. In particular,
irrational automatic numbers are transcendental.
Our main tool is a new, combinatorial transcendence criterion.
}

}

\vskip 6mm
\vskip 8mm

\centerline{\bf 1. Introduction}

\vskip 6mm

Let $b \ge 2$ be an integer. The $b$-adic expansion of every
rational number is eventually periodic, but what can be said on the
$b$-adic expansion of an irrational algebraic number? This
question was addressed for the first time by \'Emile Borel {[11]}, 
who made the conjecture that such an expansion should satisfies the same 
laws as do almost all real numbers. In particular, it is expected that every 
irrational algebraic number is normal in base $b$. Recall that a real
number $\theta$ is called {\it normal in base $b$} if, 
for any positive integer $n$,
each one of the $b^n$ blocks of length $n$
on the alphabet $\{0, 1, \ldots, b-1\}$
occurs in the $b$-adic expansion of $\theta$ with the same frequency $1/b^n$.
This conjecture is reputed to be out of reach:
we even do not know whether the digit $7$ occurs infinitely often in
the decimal expansion of $\sqrt{2}$.
However, some (very) partial
results have been established. 

As usual, we measure the complexity
of an infinite word ${\bf u} = u_1 u_2 \ldots $ 
defined on a finite alphabet by counting the
number $p(n)$ of distinct blocks of length $n$ occurring in the
word ${\bf u}$. In particular, the $b$-adic expansion of every 
real number normal in base $b$ satisfies $p(n)= b^n$ for any 
positive integer $n$.
Using a clever reformulation of a theorem of Ridout {[31]},
Ferenczi and Mauduit {[18]} 
established the transcendence of the
real numbers whose $b$-adic expansion is a
non eventually periodic sequence of minimal complexity,
that is, which satisfies $p(n) = n+1$ for every $n \ge 1$
(such a sequence is called a Sturmian sequence, 
see the seminal papers by Morse and Hedlund 
[26, 27]).
The combinatorial criterion given in {[18]} has
been used subsequently to exhibit further examples of
transcendental numbers with low complexity
{[3, 6, 4, 32]}. It also implies
that the complexity of the $b$-adic expansion of every 
irrational algebraic number 
satisfies $\liminf_{n \to \infty} (p(n) - n) = + \infty$. 
Although this is very far away from what is expected, 
no better result is known. 

In 1965, Hartmanis and Stearns {[19]} 
proposed an alternative approach
for the notion of complexity of real numbers, by emphasizing the
quantitative aspect of the notion of calculability introduced by
Turing {[39]}. According to them, a real number is
said to be computable in time $T(n)$ if there exists a multitape 
Turing machine 
which gives the first $n$-th terms of its binary 
expansion in (at most) $T(n)$ operations. The `simpler' real numbers in that
sense, that is, the numbers for which one can choose $T(n) = O(n)$,
are said to be computable in real time. Rational numbers share  
clearly this property. The problem of Hartmanis and Stearns, to
which a negative answer is expected, is the following: do there
exist irrational algebraic numbers which are computable in real time?
In 1968, Cobham {[13]} suggested to restrict this problem to a 
particular class of Turing machines, namely to the case of finite 
automata (see Section 3 for a definition). After several attempts by Cobham 
{[13]} in 1968 and 
by Loxton and van der Poorten {[21]} in 1982, 
Loxton and van der Poorten 
{[22]} finally claimed to have
completely solved the restricted problem in 1988. 
More precisely, they asserted that the $b$-adic expansion
of every irrational algebraic number cannot be generated by a finite
automaton. The proof proposed in {[22]}, which
rests on a method introduced by Mahler
{[23, 24, 25]}, 
contains unfortunately a rather serious gap, 
as explained by Becker {[8]}
(see also {[40]}). 
Furthermore, the
combinatorial criterion established in {[18]} is too weak
to imply this statement, often referred to as the conjecture of
Loxton and van der Poorten.

In the present paper, we prove new results concerning both notions 
of complexity. Our Theorem 1 provides a sharper lower estimate for the
complexity of the $b$-adic expansion of every irrational algebraic number.
We are still far away from proving that such an expansion is normal,
but we considerably improve upon the earlier known results.
We further establish (Theorem 2) the conjecture of Loxton
and van der Poorten, namely that 
irrational automatic numbers are transcendental.
Our proof yields more general statements and allows us to
confirm that irrational morphic numbers are transcendental, for a wide
class of morphisms (Theorems 3 and 4).

We derive Theorems 1 to 4 from a refinement
(Theorem 5) of the combinatorial
criterion from {[18]}, that we obtain as a 
consequence of the Schmidt Subspace Theorem.

Throughout the present paper, we adopt the following convention.
We use small letters ($a$, $u$, etc.) to denote letters from 
some finite alphabet ${\cal A}$.
We use capital letters ($U$, $V$, $W$, etc.) to denote
finite words. We use bold small letters (${\bf a}$, ${\bf u}$, etc.)
to denote infinite sequence of letters. We often identify
the sequence ${\bf a} = (a_k)_{k \ge 1}$ with the infinite word
$a_1 a_2 \ldots$, also called ${\bf a}$. This should not cause any confusion.

Our paper is organized as follows. The main results are stated in
Section 2 and proved in Section 5. Some definitions from
automata theory and combinatorics on words are recalled in Section 3.
Section 4 is devoted to
the new transcendence criterion and its proof. Finally, we show
in Section 6 that the Hensel expansion of every irrational
algebraic $p$-adic number cannot have a low complexity, and we conclude
in Section 7 by miscellaneous remarks.

Some of the results of the present paper were announced 
in {[2]}.

\vskip 3mm

\noi {\bf Acknowledgements.} We would like to thank Guy Barat and
Florian Luca for their useful comments. The first author 
is also most grateful to 
Jean-Paul Allouche and Val\'erie Berth\'e for their constant support.


\vskip 4mm

\centerline{\bf 2. Main results}

\vskip 6mm

As mentioned in the first part of the Introduction,
we measure the complexity of a real number
written in some integral base $b \ge 2$ by counting, 
for any positive integer $n$, the
number $p(n)$ of distinct blocks of $n$ digits (on the alphabet 
$\{0, 1,  \ldots, b-1\}$) occurring in its $b$-adic expansion.
The fonction $p$ is commonly called the {\it complexity function}.
It follows from results of Ferenczi and Mauduit {[18]} 
(see also {[4]}, Th\'eor\`eme 3) 
that the complexity function $p$
of every irrational algebraic number satisfies
$$
\liminf_{n\to\infty}\, (p(n) - n) =  + \infty.  \eqno (1)
$$
As far as we are aware, no better result is known, although
it has been proved {[3, 6, 32]} 
that some special real numbers with linear complexity are transcendental.

Our first result is a considerable improvement of (1).

\proclaim Theorem 1. Let $b \ge 2$ be an integer.
The complexity function of the $b$-adic expansion of 
every irrational algebraic number satisfies
$$
\liminf_{n\to\infty}\, {p(n) \over n}=+\infty.  
$$

It immediately follows from Theorem 1 that every irrational real number with
sub-linear complexity (i.e., such that
$p(n) = O(n)$) is transcendental. However, Theorem 1 is
slightly sharper, as is illustrated by an example due to 
Ferenczi {[17]}: he established the existence of a sequence
on a finite alphabet whose complexity function $p$ satisfies
$$
\liminf_{n \to \infty} \, {p(n) \over n}= 2
\quad {\rm and} \quad
\limsup_{n\to\infty}\, {p(n) \over n^t}=+\infty
\quad \hbox{for any $t>1$}.
$$

Most of the previous attempts towards a proof of 
the conjecture of Loxton and van der Poorten have been
made via the Mahler method {[21, 22, 8, 29]}. 
We stress that Becker {[8]}
established that, for any given non-eventually periodic
automatic sequence ${\bf u} = u_1 u_2 \ldots$, the real number
$\sum_{k \ge 1} \, u_k b^{-k}$ is transcendental, provided that 
the integer $b$ is sufficiently large (in terms of ${\bf u}$).
Since the complexity function $p$ of any automatic sequence
satisfies $p(n) = O(n)$ (see Cobham {[14]}), Theorem 1 confirms
straightforwardly this conjecture.

\proclaim Theorem 2. Let $b \ge 2$ be an integer.
The $b$-adic expansion of any irrational algebraic 
number cannot be generated by a finite automaton. 
In other words, irrational automatic numbers are transcendental.
 
Although Theorem 2 is a direct consequence of Theorem 1, we give in Section 5
a short proof of it, that rests on another result of Cobham {[14]}.

\medskip

Theorem 2 establishes a particular case of 
the following widely believed conjecture (see e.g. {[5]}).
The definitions of morphism, recurrent morphism, and morphic number
are recalled in Section 3.

\proclaim Conjecture. Irrational morphic numbers are transcendental.

Our method allows us to confirm this conjecture for a wide class
of morphisms.

\proclaim Theorem 3.
Binary algebraic irrational numbers cannot be generated by a 
morphism.

As observed by Allouche and Zamboni {[6]}, it follows from
{[18]} combined with a result of
Berstel and S\'e\'ebold {[9]} that binary irrational numbers
which are fixed point of a primitive morphism or of a
morphism of constant length $\ge 2$ are transcendental.
Our Theorem 3 is much more general.

Recently, by a totally different method,
Bailey, Borwein, Crandall, and Pomerance {[7]}
established new, interesting results on the density of the digits in 
the binary expansion of algebraic numbers.

For $b$-adic expansions with $b\ge 3$, we obtain 
a similar result as in Theorem 3, but  an additional assumption is needed.  

\proclaim Theorem 4.
Let $b\ge 3$ be an integer. The $b$-adic expansion of an algebraic 
irrational number cannot be generated by a recurrent morphism.

Unfortunately, we are unable to prove that ternary
algebraic numbers cannot be generated by a morphism. 
Consider for instance the fixed point 
$$
{\bf u} = 01212212221222212222212222221222 \ldots
$$
of the morphism defined by $0 \to 012$, $1 \to 12$, 
$2 \to 2$, and set $\alpha = \displaystyle\sum_{k \ge 1} \, u_k 3^{-k}$.
Our method does not apply to show the transcendence
of $\alpha$. Let us mention that this $\alpha$ is known to be transcendental:
this is a consequence of deep transcendence results 
proved in {[10]} and in {[15]}, concerning the values of theta series
at algebraic points.

\medskip

The proofs of Theorems 1 to 4 are given in Section 5.
The key point for them is a
new transcendence criterion, derived from the
Schmidt Subspace Theorem, and stated in Section 4.
Actually, we are able to deal also, under some conditions,
with non-integer bases (see Theorems 5 and 5A). 
Given a real number $\beta > 1$, we can expand in base $\beta$ every real 
number $\xi$ in $(0, 1)$ thanks to the greedy algorithm: 
we then get the $\beta$-expansion of $\xi$, 
introduced by Renyi {[30]}. 
Using Theorem 5, we easily see that
the conclusions of Theorems 1 to 4 remain true with the expansion in base $b$
replaced by the $\beta$-expansion when $\beta$ is a Pisot or a Salem number.  
Recall that a Pisot (resp. Salem) number is
a real algebraic integer $> 1$, whose complex conjugates
lie inside the open unit disc (resp. inside the closed unit disc, 
with at least one of them on the unit circle). In particular,
any integer $b \ge 2$ is a Pisot number. 
For instance, we get the 
following result.

\proclaim Theorem 1A. Let $\beta > 1$ be a Pisot or a Salem number.
The complexity function of the $\beta$-expansion of 
every algebraic number in $(0, 1)\setminus \Q(\beta)$ satisfies
$$
\liminf_{n\to\infty}\, {p(n) \over n}=+\infty.  
$$

Likewise, we can also state Theorems 2A, 3A, and 4A accordingly:
Theorems 1 to 4 deal with algebraic irrational numbers, while
Theorems 1A to 4A deal with algebraic numbers in $(0, 1)$ which do not lie
in the number field generated by $\beta$. 

Moreover,
our method also allows us to prove that $p$-adic
irrational numbers whose Hensel expansions have low complexity are
transcendental, see Section 6.

\vskip 6mm

\centerline{\bf 3. Finite automata and morphic sequences}

\vskip 6mm

In this Section, we gather classical definitions from
automata theory and combinatorics on words.

\medskip

\noindent{\it Finite automata and automatic sequences}. 
Let $k$ be an integer with $k\geq 2$. We denote by $\Sigma_k$ the set 
$\left\{0,1,\ldots,k-1\right\}$. A $k$-automaton is defined as a $6$-tuple 
$$
A=\left(Q,\Sigma_k,\delta,q_0,\Delta,\tau\right),
$$
where $Q$ is a finite set of states, $\Sigma_k$ is the input alphabet, 
$\delta:Q\times\Sigma_k\rightarrow Q$ is the transition function, $q_0$ is 
the initial state, $\Delta$ is the output alphabet and $\tau : Q\rightarrow 
\Delta$ is the output function.

For a state $q$ in $Q$ and for a finite word $W=w_1w_2\ldots w_n$
on the alphabet $\Sigma_k$, 
we define recursively $\delta(q,W)$ by 
$\delta(q,W)=\delta(\delta(q,w_1w_2\ldots w_{n-1}),w_n)$. 
Let $n\geq 0$ be an integer and let 
$w_r w_{r-1}\ldots w_1 w_0$ in $\left(\Sigma_k\right)^r$ 
be the $k$-ary expansion 
of $n$; thus, $n=\displaystyle\sum_{i=0}^r w_i k^{i}$. We denote by $W_n$ the 
word $w_0 w_1 \ldots w_r$. Then, a sequence ${\bf a}=(a_n)_{n\geq 0}$ is 
said to be $k$-automatic if there exists a $k$-automaton $A$ such that 
$a_n=\tau(\delta(q_0,W_n))$ for all $n\geq 0$.

\medskip

A classical example of a $2$-automatic sequence is given by the binary 
Thue-Morse sequence ${\bf a}=(a_n)_{n\geq 0}= 0110100110010 \ldots$
This sequence is defined 
as follows: $a_n$ is equal to $0$ (resp. to $1$)
if the sum of the digits in the 
binary expansion of $n$ is even (resp. is odd). It is easy to 
check that this sequence can be generated by the $2$-automaton
$$
A=\bigl(\{q_0, q_1\}, \{0, 1\}, \delta, q_0, \{0, 1\}, \tau \bigr),
$$
where 
$$
\delta(q_0, 0) = \delta (q_1, 1) = q_0, \qquad
\delta(q_0, 1) = \delta (q_1, 0) = q_1,
$$
and $\tau (q_0) = 0$, $\tau (q_1) = 1$.

\bigskip

\noindent{\it Morphisms}. 
For a finite set ${\cal A}$, we denote by ${\cal A}^*$ the free monoid 
generated by ${\cal A}$. The empty word $\varepsilon$ is the neutral element 
of ${\cal A}^*$. Let ${\cal A}$ and ${\cal B}$ be two finite sets. An 
application from ${\cal A}$ to ${\cal B}^*$ can be uniquely extended to an 
homomorphism between the free monoids ${\cal A}^*$ and ${\cal B}^*$. We 
call morphism from ${\cal A}$ to ${\cal B}$ such an homomorphism.

\medskip

\noindent{\it Sequences generated by a morphism}. 
A morphism 
$\phi$ from ${\cal A}$ into itself is said to be prolongable if there exists a 
letter $a$ such that $\phi(a)=aW$, where $W$ is a non-empty word such 
that $\phi^k(W)\not=\varepsilon$ for every $k\geq 0$. In that case, the 
sequence of finite words $(\phi^k(a))_{k\geq 1}$ converges in 
${\cal A}^{\N}$ (endowed with the product topology of the discrete 
topology on each copy of ${\cal A}$) 
to an infinite word ${\bf a}$. This infinite word is clearly 
a fixed point for $\phi$ and we say that ${\bf a}$ is generated by the morphism 
$\phi$. If, moreover, every letter occurring in ${\bf a}$
occurs at least twice, then we say that
${\bf a}$ is generated by a recurrent morphism. 
If the alphabet ${\cal A}$ has two letters, then we say that
${\bf a}$ is generated by a binary morphism. 
More generally, an infinite sequence ${\bf a}$ in ${\cal A}^{\N}$ 
is said to be morphic 
if there exist a sequence ${\bf u}$ generated by a morphism defined over 
an alphabet ${\cal B}$ and a morphism from ${\cal B}$ to ${\cal A}$ such that 
${\bf a}=\phi({\bf u})$.

\medskip

For instance, the Fibonacci morphism $\sigma$ defined from the alphabet 
$\{0,1\}$ into 
itself by $\sigma(0)=01$ and $\sigma(1)=0$ is 
a binary, recurrent morphism which generates the Fibonacci infinite 
word
$$
{\bf a}=\lim_{n\to\infty}\sigma^n(0)=010010100100101001\ldots
$$ 
This infinite word is an example of a Sturmian sequence and its complexity 
function satisfies 
thus $p(n)=n+1$ for every positive integer $n$. 

\bigskip

\noindent{\it Automatic and morphic real numbers}. 
Following the previous definitions, we say that a real number $\alpha$ 
is automatic 
(respectively, generated by a morphism, generated by a recurrent 
morphism, or morphic) if there exists an integer $b\geq 2$ such that the 
$b$-adic expansion of $\alpha$ is automatic (respectively, generated by a 
morphism, generated by a recurrent 
morphism, or morphic).

\medskip

A classical example of binary automatic number is given by 
$$
\sum_{n\geq 1} \, {1\over 2^{2^n}}
$$ 
which is transcendental, as proved by Kempner {[20]}.

\vskip 6mm

\centerline{\bf 4. A transcendence criterion for stammering sequences}

\vskip 6mm

First, we need to introduce some notation.
Let ${\cal A}$ be a finite set. The length of a word
$W$ on the alphabet ${\cal A}$, that is, the number of letters
composing $W$, is denoted by $\vert W\vert$.
For any positive integer $\ell$, we write
$W^{\ell}$ for the word $W\ldots W$ ($\ell$ times repeated concatenation
of the word $W$). More generally, for any positive real number
$x$, we denote by $W^x$ the word
$W^{\lfloor x \rfloor}W'$, where $W'$ is the prefix of
$W$ of length $\left\lceil(x-\lfloor x\rfloor)\vert W\vert\right\rceil$. 
Here, and in all what follows, $\lfloor y\rfloor$ and
$\lceil y\rceil$ denote, respectively, the integer part and the upper
integer part of the real number $y$. 
Let ${\bf a}=(a_k)_{k \ge 1}$ be a sequence of elements from ${\cal A}$,
that we identify with the infinite word $a_1 a_2 \ldots$
Let $w>1$ be a real number.
We say that ${\bf a}$ 
satisfies Condition $(*)_w$ if ${\bf a}$ is not
eventually periodic and if there exist 
two sequences of finite words $(U_n)_{n \ge 1}$,   
$(V_n)_{n \ge 1}$ such that:

\medskip

\item{\rm (i)} For any $n \ge 1$, the word $U_nV_n^w$ is a prefix
of the word ${\bf a}$;

\smallskip

\item{\rm (ii)} The sequence
$({\vert U_n\vert} / {\vert V_n\vert})_{n \ge 1}$ is bounded from above;

\smallskip

\item{\rm (iii)} The sequence $(\vert V_n\vert)_{n \ge 1}$ is 
increasing.

\medskip

As suggested to us by Guy Barat, a sequence satisfying Condition $(*)_w$
for some $w>1$ may be called a stammering sequence.

\proclaim Theorem 5. Let $\beta > 1$ be a Pisot or a Salem number.
Let ${\bf a}=(a_k)_{k \ge 1}$ 
be a bounded sequence of rational integers. 
If there exists a real number $w>1$
such that ${\bf a}$ satisfies Condition $(*)_w$, then the real
number 
$$
\alpha:= \sum_{k =1}^{+\infty}\, {a_k \over \beta^k}
$$ 
either belongs to $\Q(\beta)$, or is transcendental.

The proof of Theorem 5 rests on the Schmidt Subspace Theorem 
{[37]} (see also {[38]}), and more precisely
on a $p$-adic generalization due to Schlickewei {[34, 35]}
and Evertse {[16]}. Note that the particular case when $\beta$ is an
integer $\ge 2$ was proved in {[2]}.
Note also that Adamczewski [1] proved that, under a stronger
assumption on the sequence $(a_k)_{k \ge 1}$, the number $\alpha$
defined in the statement of Theorem 5 is transcendental.

\bigskip

\noi {\bf Remarks}.  

${\bullet}$ Theorem 5 is considerably stronger than the 
criterion of Ferenczi and Mauduit {[18]}: 
our assumption $w > 1$ replaces their
assumption $w > 2$. This type of condition is rather flexible,
compared with the Mahler method, for which a functional equation 
is needed. For instance, the conclusion of Theorem 5 also
holds if the sequence ${\bf a}$ is an unbounded 
sequence of integers that does not increase
too rapidly. Nevertheless, one should acknowledge that,
when it can be applied, the Mahler method gives the transcendence
of the infinite series $\sum_{k=1}^{+\infty}\, a_k \beta^{-k}$
for every algebraic number $\beta$ such that this series converges.

\medskip

${\bullet}$ We emphasize that if a sequence ${\bf u}$ 
satisfies Condition $(*)_w$ and if $\phi$ is a non-erasing morphism
(that is, if the image by $\phi$ of any letter has length at least $1$),
then $\phi({\bf u})$ satisfies Condition $(*)_w$, as well.
This observation is used in the proof of Theorem 2.

\medskip

${\bullet}$ If $\beta$ is an algebraic number which is neither a
Pisot, nor a Salem number, it is still possible to get
a transcendence criterion using the approach followed for
proving Theorem 5. However, the
assumption $w>1$ should then 
be replaced by a weaker one, involving the Mahler measure of $\beta$
and $\limsup_{n\to\infty} \, {\vert U_n\vert} / {\vert V_n\vert}$.
Furthermore, the same approach shows
that the full strength of Theorem 5 holds when $\beta$ is a Gaussian
integer. More details will be given in a subsequent work.

\bigskip

Before beginning the proof of Theorem 5, we quote a version of the
Schmidt Subspace Theorem, as formulated by Evertse {[16]}.

We normalize absolute values and heights as follows. 
Let $\K$ be an algebraic number field of degree $d$. Let $M(\K)$ denote 
the set of places on $\K$.
For $x$ in $\K$ and a place $v$ in $M(\K)$, define the
absolute value $|x|_v$ by
\medskip

\item{\rm (i)} $|x|_v = |\sigma (x)|^{1/d}$ \quad \hbox{if
$v$ corresponds to the embedding $\sigma : \K \hookrightarrow \R$};

\smallskip
\item{\rm (ii)} $|x|_v = |\sigma (x)|^{2/d} = 
|\overline{\sigma} (x)|^{2/d}$ \quad \hbox{if
$v$ corresponds to the pair of conjugate complex em-}
\hbox{beddings 
$\sigma , \overline{\sigma} : \K \hookrightarrow \C$};

\smallskip
\item{\rm (iii)} $|x|_v = (N \pp)^{-ord_\pp (x) / d}$ \quad \hbox{if
$v$ corresponds to the prime ideal $\pp$ of $O_{\K}$}.

These absolute values satisfy the product formula
$$
\prod_{v \in M(\K)} \, |x|_v = 1 \qquad
\hbox{for $x$ in $\K^*$}.
$$

Let ${\bf x} = (x_1, \ldots, x_n)$ be in
$\K^n$ with ${\bf x} \not= {\bf 0}$.
For a place $v$ in $M(\K)$, put
$$
\eqalign{
|{\bf x}|_v & = \biggl( \, \sum_{i=1}^n \, |x_i|_v^{2d}
\, \biggr)^{1/ (2d)} \quad \hbox{if $v$ is real infinite}; \cr
|{\bf x}|_v & = \biggl( \, \sum_{i=1}^n \, |x_i|_v^{d}
\, \biggr)^{1/ d} \quad \hbox{if $v$ is complex infinite}; \cr
|{\bf x}|_v & = \max\{ |x_1|_v, \ldots, |x_n|_v \}
\quad \hbox{if $v$ is finite}. \cr}
$$
Now define the {\it height} of ${\bf x}$ by
$$
H({\bf x}) = H(x_1, \ldots, x_n) = \prod_{v \in M(\K)} \, |{\bf x}|_v.
$$
We stress that $H({\bf x})$ depends only on ${\bf x}$ and not
on the choice of the number field $\K$ containing 
the coordinates of ${\bf x}$, see e.g. {[16]}.

We use the following formulation of the Subspace Theorem over number fields.
In the sequel, we assume that the algebraic closure of $\K$
is ${\overline{\bf Q}}$. We choose for every place $v$ in $M(\K)$ 
a continuation of $| \cdot |_v$ to ${\overline{\bf Q}}$, that we denote
also by $| \cdot |_v$.

\proclaim Theorem E. Let $\K$ be an algebraic number field. 
Let $m\ge 2$ be an integer.
Let $S$ be a finite set of places on $\K$ containing all
infinite places. For each $v$ in $S$, let $L_{1, v}, \ldots , L_{m, v}$ 
be linear forms with algebraic coefficients and with
$$
{\it rank} \, \{L_{1, v}, \ldots , L_{m, v}\} = m.
$$
Let $\eps$ be real with $0 < \eps < 1$.
Then, the set of solutions ${\bf x}$
in $\K^m$ to the inequality
$$
\prod_{v \in S} \, \prod_{i=1}^m \,
{\vert L_{i, v} ({\bf x}) \vert_v \over 
\vert {\bf x} \vert_v} \le \, H({\bf x})^{-m-\eps}
$$
lies in finitely many proper subspaces of $\K^m$.

For a proof of Theorem E, the reader is directed to {[16]},
where a quantitative version is established (in the sense that one
bounds explicitly the number of exceptional subspaces).

\medskip

We now turn to the proof of Theorem 5.
Keep the notation and the assumptions of this theorem.
Assume that the parameter $w > 1$ is fixed, as well as the 
sequences $(U_n)_{n \ge 1}$ and $(V_n)_{n \ge 1}$ occurring in the
definition of Condition $(*)_w$. 
Set also $r_n=\vert U_n\vert$ and $s_n=\vert V_n\vert$ for any $n \ge 1$.
We aim to prove that the real number
$$
\alpha:= \sum_{k=1}^{+\infty}\, {a_k \over \beta^k}
$$ 
either lies in $\Q(\beta)$ or is transcendental. 
The key fact is the observation
that $\alpha$ admits infinitely many good approximants in the number
field $\Q(\beta)$ 
obtained by truncating its expansion and completing it by periodicity.
Precisely, for any positive integer $n$, we define the sequence
$(b_k^{(n)})_{k \ge 1}$ by
$$
\eqalign{
b_h^{(n)} & = a_h \quad \hbox{for $1 \le h \le r_n + s_n$,} \cr
b_{r_n + h + j s_n}^{(n)} & = a_{r_n + h} 
\quad \hbox{for $1 \le h \le s_n$ and $j \ge 0$.} \cr}
$$
The sequence
$(b_k^{(n)})_{k \ge 1}$ is eventually periodic, with preperiod $U_n$
and with period $V_n$. Set
$$
\alpha_n= \sum_{k=1}^{+\infty}\, {b_k^{(n)} \over \beta^k},
$$
and observe that
$$
\alpha - \alpha_n = \sum_{k=r_n + \lceil w s_n \rceil +1}^{+\infty}
\, {a_k - b_k^{(n)} \over \beta^k}\cdot \eqno (2)
$$

\proclaim Lemma 1. 
For any integer $n$, there exists an integer polynomial $P_n(X)$
of degree at most $r_n+s_n-1$ such that
$$
\alpha_n=\, {P_n(\beta) \over \beta^{r_n}(\beta^{s_n}-1)}\cdot
$$
Further, the coefficients of $P_n (X)$
are bounded in absolute value by $2 \max_{k \ge 1} |a_k|$.

\medskip

\noindent {\bf proof.} By definition of $\alpha_n$, we get
$$
\eqalign{
\alpha_n & = \sum_{k=1}^{r_n}\, {a_k \over \beta^k}+
\sum_{k=r_n + 1}^{+\infty}\, {b_k^{(n)} \over \beta^k} 
=\sum_{k=1}^{r_n}\, {a_k \over \beta^k}+ {1 \over \beta^{r_n}} \, 
\sum_{k=1}^{+\infty}\, {b_{r_n+k}^{(n)} \over \beta^k} \cr
& = \sum_{k=1}^{r_n}\, {a_k \over \beta^k}+ {1 \over \beta^{r_n}}
\sum_{k=1}^{s_n} \, {a_{r_n+k } \over \beta^k} \, \biggl(
\sum_{j=0}^{+\infty} \, {1 \over \beta^{j s_n}} \biggr) \cr
& = \sum_{k=1}^{r_n}\, {a_k \over \beta^k}+
\sum_{k=1}^{s_n}\, {a_{r_n+k } \over \beta^{k + r_n  - s_n}
(\beta^{s_n}-1)} = {P_n (\beta) \over \beta^{r_n}(\beta^{s_n}-1)}, \cr}
$$
where we have set
$$
P_n (X) = \sum_{k=1}^{r_n} \, a_k \, X^{r_n-k}(X^{s_n}-1) +
\sum_{k=1}^{s_n} \, a_{r_n+k } \,  X^{s_n - k}. 
$$
The last assertion of the lemma is clear.  \cqfd

Set $\K = \Q(\beta)$ and denote by $d$ the degree of $\K$.
We assume that $\alpha$ is algebraic, and we consider the following
linear forms, in three variables and with algebraic coefficients.
For the place $v$ corresponding to the embedding of $\K$ defined by
$\beta \hookrightarrow \beta$, set
$L_{1, v}(x,y,z)= x$, $L_{2, v}(x,y,z)=y$, and
$L_{3, v}(x,y,z)=\alpha x+\alpha y+z$. 
It follows from (2) and Lemma 1 that
$$
\vert L_{3, v}(\beta^{r_n+s_n},-\beta^{r_n},- P_n(\beta))\vert_v   
= \vert \alpha (\beta^{r_n}(\beta^{s_n}-1))- P_n(\beta) \vert^{1/d}
\ll {1 \over \beta^{(w-1)s_n/d}},  \eqno (3)
$$
where we have chosen
the continuation of $| \cdot |_v$ to ${\overline{\bf Q}}$
defined by $| x |_v = |x|^{1/d}$.
Here and throughout this Section, the constants implied by the
Vinogradov symbol $\ll$ depend (at most) on $\alpha$, $\beta$, and
$\max_{k \ge 1} |a_k|$, but are independent of $n$.

Denote by $S'_{\infty}$ the set of all other infinite places
on $\K$ and by $S_0$ the set of all finite places on $\K$ dividing $\beta$.  
Observe that $S_0$ is empty if $\beta$ is an algebraic unit. 
For any $v$ in $S_0 \cup S'_{\infty}$, set
$L_{1, v}(x,y,z)=x$, $L_{2, v}(x,y,z)=y$, and $L_{3, v}(x,y,z)=z$. 
Denote by $S$ the union of $S_0$ and the infinite places on $\K$.
Clearly, for any $v$ in $S$, the forms $L_{1, v}$, $L_{2, v}$ and $L_{3, v}$
are linearly independent. 

To simplify the exposition, set
$$
{\bf x}_n = (\beta^{r_n+s_n},-\beta^{r_n},- P_n(\beta)).
$$
We wish to estimate the product
$$
\Pi := \prod_{v \in S} \, \prod_{i=1}^3 \,
{\vert L_{i, v} ({\bf x}_n) \vert_v \over 
\vert {\bf x}_n \vert_v} = \prod_{v \in S}  \, |\beta^{r_n + s_n}|_v
\, |\beta^{r_n}|_v \, {|L_{3, v} ({\bf x}_n)|_v \over |{\bf x}_n|_v^3}
$$
from above. By the product formula and the definition of $S$, we 
immediately get that
$$
\Pi = \prod_{v \in S} \, {|L_{3, v} ({\bf x}_n)|_v \over |{\bf x}_n|_v^3}.
\eqno (4)
$$
Since the polynomial $P_n (X)$ has integer coefficients and since $\beta$ is
an algebraic integer, we have
$|L_{3, v} ({\bf x}_n)|_v = |P_n (\beta)|_v \le 1$ 
for any place $v$ in $S_0$.
Furthermore, as the conjugates of $\beta$ have moduli at most $1$,
we have for any infinite place $v$ in $S'_{\infty}$
$$
|L_{3, v} ({\bf x}_n)|_v \ll (r_n + s_n)^{d_v/d},
$$
where $d_v = 1$ or $2$ according as $v$ is real infinite or complex infinite,
respectively. Together with (3) and (4), this gives
$$
\eqalign{
\Pi & \ll (r_n + s_n)^{(d-1) / d} \, \beta^{-(w-1)s_n/d} \,
\prod_{v \in S} \, |{\bf x}_n|_v^{-3} \cr
& \ll (r_n + s_n)^{(d-1) / d} \, \beta^{- (w-1)s_n/d} \, 
H ( {\bf x}_n )^{-3}, \cr}
$$
since $|{\bf x}_n|_v = 1$ if $v$ does not belong to $S$.

Furthermore, it follows from Lemma 1 and from the fact that the moduli of the
complex conjugates of $\beta$ are at most $1$ that
$$
H ( {\bf x}_n ) \ll (r_n + s_n)^d \, \beta^{(r_n + s_n)/d}.
$$
Consequently, we infer from Condition $(*)_w$ that
$$
\eqalign{
\prod_{v \in S} \, \prod_{i=1}^3 \,
{\vert L_{i, v} ({\bf x}_n) \vert_v \over 
\vert {\bf x}_n \vert_v} & \ll (r_n + s_n)^{d w} \, 
H ( {\bf x}_n )^{- (w-1) s_n / (r_n + s_n)} \, 
H ( {\bf x}_n )^{-3} \cr & \ll H ( {\bf x}_n )^{-3 - \eps}, \cr}
$$
for some positive real number $\eps$. 

It then follows from
Theorem E that the points $(\beta^{r_n+s_n},-\beta^{r_n},- P_n(\beta))$
lie in a finite number of proper subspaces of $\K^3$.
Thus, there exist a non-zero triple $(x_0,y_0,z_0)$ in $\K^3$ and
infinitely many integers $n$ such that
$$
x_0-y_0 {\beta^{r_n} \over \beta^{r_n+s_n}}-
z_0 {P_n(\beta) \over \beta^{r_n+s_n}}=0.
$$
Taking the limit along this subsequence of
integers and noting that $(s_n)_{n \ge 1}$ 
tends to infinity, we get that $x_0 = z_0 \alpha$.
Thus, $\alpha$ belongs to $\K=\Q (\beta)$, as asserted. \cqfd

\medskip

Let us restrict our attention to the case when $\beta$
is a Pisot number. 
Dealing with the $\beta$-expansions of real numbers 
(instead of arbitrary power series in $\beta$) allows us to 
improve the conclusion of Theorem 5. 

\proclaim Theorem 5A.
Let $\beta > 1$ be a Pisot number. 
Let $\alpha$ be in $(0, 1)$, and
consider its $\beta$-expansion
$$
\alpha:= \sum_{k =1}^{+\infty}\, {a_k \over \beta^k}.
$$ 
If $(a_k)_{k \ge 1}$ satisfies Condition $(*)_w$ for some
real number $w>1$, then $\alpha$ is transcendental.

\medskip

\noindent{\bf Proof.} By a result of K. Schmidt {[36]},
we know that the $\beta$-expansion of every element of
$\Q(\beta) \cap (0, 1)$ is
eventually periodic. Thus, it does not
satisfy Condition $(*)_w$. We conclude by applying
Theorem 5. \cqfd

Note that for a Salem number $\beta$,
it is an important open problem to decide whether every element of
$\Q(\beta) \cap (0, 1)$ has an
eventually periodic $\beta$-expansion.


\vskip 6mm

\centerline{\bf 5. Proofs of Theorems 1 to 4}

\vskip 6mm

We begin by a short proof of Theorem 2.

\bigskip

\noi {\bf Proof of Theorem 2}. 
Let ${\bf a} = (a_k)_{k \ge 1}$ be a non-eventually periodic 
automatic sequence defined on a finite alphabet ${\cal A}$. 
Recall that a morphism is
called uniform if the images of each letter have the same length. 
Following Cobham {[14]}, 
there exist a morphism $\phi$ from 
an alphabet ${\cal B}=\{1,2,\ldots,r\}$ 
to the alphabet ${\cal A}$ and an uniform
morphism $\sigma$ from ${\cal B}$ into itself such that ${\bf
a}=\phi({\bf u})$, where ${\bf u}$ is a fixed point for ${\sigma}$. 
Observe first that the sequence ${\bf a}$ satisfies 
Condition $(*)_w$ if this is the case for ${\bf u}$. 
Further, by the Dirichlet {\it Schubfachprinzip}, the prefix of length
$r+1$ of ${\bf u}$ can be written under the form
$W_1 u W_2 u W_3$, where $u$ is a letter
and $W_1$, $W_2$, $W_3$ are (possibly empty) finite words.
We check that the assumptions of Theorem 1 are satisfied
by ${\bf u}$ with the sequences $(U_n)_{n \ge 1}$ and $(V_n)_{n \ge 1}$
defined for any $n \ge 1$ by $U_n = \sigma^n (W_1)$ 
and $V_n = \sigma^n (u W_2)$. Indeed, since $\sigma$ is 
a morphism of constant length, we get, on the one hand, that
$$
{|U_n| \over |V_n|} \le {|W_1| \over 1 + |W_2|} \le r-1
$$
and, on the other hand, that
$\sigma^n (u)$ is a prefix of $V_n$ of length at least
$1/r$ times the length of $V_n$. 
It follows that Condition $(*)_{1+1/r}$ is satisfied
by the sequence ${\bf u}$, and thus by our sequence ${\bf a}$
(here, we use the observation we made in Section 4). 
Let $b \ge 2$ be an integer. By applying Theorem 5 with $\beta = b$,
we conclude that the automatic number $\sum_{k = 1}^{+ \infty} \,
a_k b^{-k}$ is transcendental.  \cqfd

\bigskip

\noi {\bf Proof of Theorem 1}. 
Let $\alpha$ be an irrational number. Without any loss of
generality, we assume that $\alpha$ is in $(0, 1)$
and we denote by $0.u_1u_2\ldots u_k \ldots$ its $b$-adic
expansion. The sequence $(u_k)_{k \ge 1}$ takes its values in
$\{0,1,\ldots,b-1\}$ and is not eventually periodic. We assume
that there exists an integer $\kappa \ge 2$ such that the complexity
function $p$ of $(u_k)_{k \ge 1}$ satisfies
$$
p(n) \le \kappa n \qquad \hbox{for infinitely many integers $n \ge 1$},
$$
and we shall derive that Condition $(*)_w$
is then fulfilled by the sequence $(u_k)_{k \ge 1}$
for a suitable $w>1$. By
Theorem 5, this will imply that $\alpha$ is transcendental.

Let $n_k$ be an integer with $p(n_k) \le \kappa n_k$.
Denote by $U(\ell)$ the prefix of 
${\bf u} := u_1 u_2 \ldots $ of length $\ell$.
By the Dirichlet {\it Schubfachprinzip}, there exists (at least)
one word $M_k$ of length $n_k$
which has (at least) two occurrences in $U ((\kappa+1) n_k)$.
Thus, there are (possibly empty) words $A_k$, $B_k$, $C_k$ and $D_k$, 
such that
$$
U ((\kappa+1) n_k) = A_k M_k C_k D_k = A_k B_k M_k D_k \qquad
{\rm and} \qquad |B_k| \ge 1.
$$
We observe that $|A_k| \le \kappa n_k$.
We have to distinguish three cases:

\item{\rm (i)} $|B_k| > |M_k|$;

\smallskip
\item{\rm (ii)} $\lceil |M_k|/3 \rceil \le |B_k| \le |M_k|$;

\smallskip
\item{\rm (iii)} $1 \le |B_k| < \lceil |M_k|/3 \rceil$.

\bigskip
${\rm (i)}$. Under this assumption, there exists a word $E_k$ such that
$$
U ((\kappa+1) n_k) = A_k M_k E_k M_k D_k.
$$
Since $|E_k| \le (\kappa-1) |M_k|$, the word $A_k (M_k E_k)^s$ with
$s = 1 + 1/\kappa$ is a prefix of ${\bf u}$.
Furthermore, we observe that
$$
|M_k E_k| \ge |M_k| \ge {|A_k| \over \kappa }.
$$

\bigskip

${\rm (ii)}$. Under this assumption, there exist two words $E_k$ and $F_k$ 
such that
$$
U ((\kappa+1) n_k) = A_k M_k^{1 / 3} E_k 
M_k^{1 / 3} E_kF_k.
$$
Thus, the word 
$A_k (M_k^{1 / 3} E_k)^2$ is a prefix of ${\bf u}$.
Furthermore, we observe that
$$
|M_k^{1 / 3} E_k| \ge {|M_k| \over 3}
\ge {|A_k| \over 3\kappa }.
$$

\bigskip

${\rm (iii)}$. In the present case, $B_k$ is clearly a prefix of $M_k$, and we infer
from $B_k M_k = M_k C_k$ that $B_k^t$ is a 
prefix of $M_k$, where $t$ is the integer part of $|M_k|/|B_k|$. 
Observe that $t \ge 3$.
Setting $s = \lfloor t / 2\rfloor$, we see that
$A_k (B_k^s)^2$ is a prefix of ${\bf u}$ and
$$
|B_k^s| \ge {|M_k| \over 4} \ge {|A_k| \over 4\kappa }.
$$

\medskip

In each of the three cases above, we have proved that there are
finite words $U_k$, $V_k$ such that $U_k V_k^{1 + 1/\kappa}$ 
is a prefix of ${\bf u}$ and:  

\medskip

\item{$\bullet$} $\vert U_k\vert\leq \kappa n_k$;

\smallskip
\item{$\bullet$} $\vert V_k\vert\geq n_k / 4$;

\smallskip
\item{$\bullet$} $w \geq 1+ 1/\kappa >1$.

\medskip
\noindent 
Consequently, the sequence 
$({\vert U_k\vert / \vert V_k\vert} )_{k\geq 1}$ 
is bounded from above by $4\kappa$. Furthermore, it follows from
the lower bound $\vert V_k\vert\geq n_k/4$ that
we may assume that the sequence $(\vert V_k\vert)_{k\geq 1}$
is strictly increasing. This implies that the sequence ${\bf u}$ 
satisfies Condition $(*)_{1 + 1/\kappa}$. By applying
Theorem 5 with $\beta = b$,
we conclude that $\alpha$ is transcendental. \cqfd

\bigskip

\noi {\bf Proof of Theorem 3}. 
Let ${\bf a}$ be a sequence generated by a morphism
$\phi$ defined on a finite alphabet ${\cal A}$. For any positive integer
$n$, there exists a letter $a_n$ satisfying 
$$
\vert\phi^n(a_n)\vert=\max\{\vert\phi^n(j)\vert : j\in{\cal A}\}.
$$  
This implies the existence of a letter $a$ in ${\cal A}$ and 
of a strictly increasing sequence
of positive integers $(n_k)_{k\geq 1}$ such that for every $k \ge 1$ we have
$$
\vert\phi^{n_k}(a)\vert=\max\{\vert\phi^{n_k}(j)\vert : j\in{\cal A}\}.
$$ 
Assume from now on that ${\cal A}$ has two elements.
Since the sequence ${\bf a}$ is not eventually periodic there exist 
at least two occurrences in ${\bf a}$ of the two elements of ${\cal A}$. 
In particular, there exist at least two occurrences of the letter $a$ in the 
sequence ${\bf a}$.  
We can thus find two (possibly empty) finite words $W_1$
and $W_2$ such that $W_1aW_2a$ is a prefix of ${\bf a}$. 
We check that the assumptions of Theorem 5 are satisfied
by ${\bf a}$ with the sequences $(U_k)_{k \ge 1}$ and $(V_k)_{k \ge 1}$
defined by $U_k = \phi^{n_k} (W_1)$ and $V_k = \phi^{n_k} (a W_2)$ 
for any $k \ge 1$. Indeed, by definition of $a$, we have
$$
{|U_k| \over |V_k|} \le {|W_1|} 
$$
and $\phi^{n_k} (a)$ is a prefix of $V_k$ of length at least
$1/(\vert W_2\vert+1)$ times the length of $V_k$. It follows that Condition 
$(*)_w$ is satisfied by the sequence ${\bf a}$ with $w=1+1/(\vert
W_2\vert+1)$. We conclude by applying Theorem 5.
\cqfd

\bigskip

\noi {\bf Proof of Theorem 4}. 
Let ${\bf a}$ be a sequence generated by a recurrent morphism
$\phi$ defined on an alphabet ${\cal A}$. 
As we have already noticed in the beginning of
the proof of Theorem 3, there exist a 
letter $a$ and a strictly increasing sequence
of positive integers $(n_k)_{k\geq 1}$ such that for every $k\ge 1$ we have
$$
\vert\phi^{n_k}(a)\vert=\max\{\vert\phi^{n_k}(j)\vert : j\in{\cal A}\}.
$$ 
Since by assumption 
the sequence ${\bf a}$ is recurrent there exist at least two 
occurrences of the letter $a$. We then apply the same trick as in the 
proof of Theorem 3, and we again conclude by applying Theorem 5.
\cqfd

 
\vskip 6mm

\goodbreak

\centerline{\bf 6. Transcendence of $p$-adic numbers}

\vskip 6mm

Let $p$ be a prime number. As usual, we denote by $\Q_p$ the field of
$p$-adic numbers. We call algebraic (resp. transcendental) any element of
$\Q_p$ which is algebraic (resp. transcendental) over $\Q$. 
A suitable version of the Schmidt Subspace Theorem,
due to Schlickewei {[33]}, 
can be applied to derive a lower bound for the complexity 
of the Hensel expansion of every irrational algebraic number in $\Q_p$.

\proclaim Theorem 1B. Let $\alpha$ be an irrational algebraic 
number in $\Q_p$
and denote by $$\alpha = \sum_{k=-m}^{+\infty} \, a_k p^k$$ its Hensel
expansion. Then, the complexity function $p$ of the sequence 
$(a_k)_{k \ge -m}$ satisfies
$$
\liminf_{n\to\infty}\, {p(n) \over n} = +\infty.
$$

Likewise (see Section 2), we can also state Theorems 2B, 3B and  4B. 
Theorem 1B follows from Theorem 6 below, along
with the arguments used in the proof of Theorem 1.

The method of proof of Theorem 5 applies to provide us with a new 
transcendence criterion for $p$-adic numbers.

\proclaim Theorem 6. Let $p$ be a prime number and let
$(a_k)_{k \ge -m}$ be a sequence taking its values in 
$\{0,1,\ldots,p-1\}$. Let $w>1$ be a real number.
If the sequence $(a_k)_{k \ge 1}$ 
satisfies Condition $(*)_w$, then the $p$-adic number
$$
\alpha:= \sum_{k=-m}^{+\infty} \, a_k p^k
$$
is transcendental.

We briefly outline the proof of Theorem 6.
Let $p$ and $(a_k)_{k \ge -m}$ be as in the statement of 
this theorem.
There exist a parameter $w>1$ and two sequences $(U_n)_{n \ge 1}$ 
and $(V_n)_{n \ge 1}$ of finite words 
as in the definition of Condition $(*)_w$. For any $n \ge 1$,
set $r_n=\vert U_n\vert$ and $s_n=\vert V_n\vert$.
To establish Theorem 6, it is enough to prove that the $p$-adic number
$$
\alpha':= \sum_{k=1}^{+\infty}\, a_k p^k
$$ 
is transcendental. As in the proof of Theorem 5,
the key fact is the observation
that $\alpha'$ admits infinitely many good rational approximants
obtained by truncating its Hensel expansion and completing by periodicity.
Precisely, for any positive integer $n$, we define the sequence
$(b_k^{(n)})_{k \ge 1}$ exactly as in Section 4, and we set
$$
\alpha_n= \sum_{k=1}^{+\infty}\, b_k^{(n)} p^k.
$$
An easy calculation shows that we have
$$
|\alpha' - \alpha_n|_p \le p^{-r_n - w s_n},
\quad \alpha_n= {p_n \over p^{s_n}-1},
$$
where
$$
p_n= \biggl(\sum_{k=1}^{r_n} a_k \, p^k \biggr)(p^{s_n}-1)-
\sum_{k=1}^{s_n} a_{r_n+k} \, p^{r_n + k}.
$$

Assuming that $\alpha'$ is an algebraic number in $\Q_p$,
we apply Theorem 4.1 of Schlickewei {[33]} with the linear forms
$L_{1, p}(x,y,z)= x$, $L_{2, p}(x,y,z)=y$, 
$L_{3, p}(x,y,z)=\alpha' x+\alpha' y+z$, 
$L_{1, \infty}(x,y,z)= x$, $L_{2, \infty}(x,y,z)=y$, and
$L_{3, \infty}(x,y,z)= z$. Setting ${\bf x}_n := (p^{s_n}, - 1, -p_n)$,
we get $|L_{1, p}({\bf x}_n)|_p = p^{-s_n}$ and
$|L_{3, p}({\bf x}_n)|_p \le p^{-r_n - w s_n}$. 
We then follow the same lines as in the proof of Theorem 5,
and we end up in a contradiction. This proves that $\alpha'$
is transcendental. \cqfd

 
\vskip 6mm

\centerline{\bf 7. Concluding remarks}

\vskip 6mm

It is of interest to compare our result with a celebrated theorem
of Christol, Kamae, Mend\`es France, and Rauzy {[12]} 
concerning algebraic elements of the field ${\bf F}_p ((X))$. 
Their result asserts that, for any given prime number $p$, 
the sequence of integers ${\bf u}=(u_k)_{k\ge 1}$ is
$p$-automatic if and only if the formal power series $\displaystyle
\sum_{k \ge 1} u_k X^k$
is algebraic over the field of rational functions ${\bf F}_p (X)$. 
Thanks to Theorem 2, we thus easily derive the following statement. 

\proclaim Theorem 7.  Let $b\ge 2$ be an integer and $p$ be a prime number. 
The formal power series $\displaystyle
\sum_{k \ge 1} u_k X^k$ and 
the real number $\displaystyle
\sum_{k \ge 1} {u_k\over b^k}$ are both algebraic 
(over ${\bf F}_p (X)$ and over $\Q$, respectively) 
if and only if they are rational.

Note that Theorem 2B (see Section 6) naturally gives 
rise to a similar result where the real number 
$\displaystyle\sum_{k \ge 1} {u_k\over b^k}$ is 
replaced by the $q$-adic number 
$\displaystyle\sum_{k \ge 1} u_kq^k$, for an arbitrary prime number $q$.
In particular, this holds true for $q=p$. 

\bigskip

In 1991, Morton and Mourant {[28]} proved the following result: 
If $k\geq 2$ is an integer, $P$ is 
a non-zero pattern of digits in base $k$, and if
$e_{k,P,b}(n)\in\{0,1,\ldots,b-1\}$ counts the number of occurrences 
modulo $b$ of $P$ in the $k$-ary expansion of $n$, then the real number 
$\alpha(k,P,b)=\displaystyle\sum_{n=0}^{+\infty}{e_{k,P,b}(n)\over b^n}$ is  
transcendental except when $k=3$, $P=1$ and $b=2$. Moreover, 
we have $\alpha(3,1,2) = 2/3$ in this particular case. 

The proof given by Morton and Mourant is based on Theorem 2 and their paper 
refers to the work of Loxton and van der Poorten {[22]}. 
The present work validates their result.

It is interesting to remark that the simplest case $k=2$, $P=1$ and $b=2$ 
corresponds to the well-known Thue--Morse number, whose transcendence has 
been proved by Mahler {[23]}. 
The theorem of Morton and Mourant can thus be seen as a 
full generalisation of the Mahler result.

\vskip 9mm

\centerline{\bf References}

\vskip 7mm

\sm \item{[1]}
B. Adamczewski, 
{\it Transcendance \og \`a la Liouville \fg de certain nombres r\'eels},
C. R. Acad. Sci. Paris 338 (2004), 511--514.

\sm \item{[2]}
B. Adamczewski, Y. Bugeaud \&  F. Luca,
{\it Sur la complexit\'e des nombres alg\'ebriques},
C. R. Acad. Sci. Paris 339 (2004), 11--14.

\sm \item{[3]} 
B. Adamczewski \& J. Cassaigne,
{\it On the transcendence of real numbers with a regular expansion},
J. Number Theory 103 (2003), 27--37.

\sm \item{[4]}
J.-P. Allouche,
{\it Nouveaux r\'esultats de transcendance de r\'eels \`a 
d\'eveloppements non al\'eatoire},
Gaz. Math. 84 (2000), 19--34.

\sm \item{[5]}
J.-P. Allouche \& J. Shallit, 
Automatic Sequences: Theory, Applications, Generalizations, 
Cambridge University Press, Cambridge, 2003.

\sm \item{[6]}
 J.-P. Allouche \& L. Q. Zamboni,
 {\it Algebraic irrational binary numbers cannot be fixed points of
              non-trivial constant length or primitive morphisms},
 J. Number Theory 69 (1998), 119--124.

\sm \item{[7]}
D. H. Bailey, J. M. Borwein, R. E. Crandall \& C. Pomerance,
{\it On the binary expansions of algebraic numbers},
J. Th\'eor. Nombres Bordeaux 16 (2004), 487--518.

\sm \item{[8]}
P. G. Becker, 
{\it $k$-regular power series and Mahler-type functional equations}, 
J. Number Theory 49 (1994), 269--286.

\sm \item{[9]}
J. Berstel \& P. S\'e\'ebold,
{\it A characterization of overlap-free morphisms},
Disc. Appl. Math. 46 (1993), 275--281.

\sm \item{[10]}
D. Bertrand,
{\it Theta functions and transcendence},
Ramanujan J. 1 (1997), 339--350.

\sm \item{[11]}
   \'E. Borel,
 {\it Sur les chiffres d\'ecimaux de $\sqrt{2}$ et divers
              probl\`emes de probabilit\'es en cha\^\i ne},
 C.~ R.~ Acad.~ Sci.~ Paris 230 (1950), 591--593.

\sm \item{[12]}
G. Christol, T. Kamae, M. Mend\`es France \& G. Rauzy,
{\it Suites alg\'ebriques, automates et substitutions},
Bull. Soc. Math. France 108 (1980), 401--419.

\sm \item{[13]}
 A. Cobham,
 {\it On the Hartmanis-Stearns problem for a class of tag machines},
 Conference Record of 1968 Ninth Annual Symposium on Switching and 
Automata Theory, Schenectady, New York (1968), 51--60.

\sm \item{[14]}
 A. Cobham,
 {\it Uniform tag sequences},
Math. Systems Theory 6 (1972), 164--192.

\sm \item{[15]}
D. Duverney, Ke. Nishioka, Ku. Nishioka and I. Shiokawa,
{\it Transcendence of Jacobi's theta series},
Proc. Japan Acad. Ser. A 72 (1996), 202--203.

\sm \item{[16]}
J.-H. Evertse,
{\it An improvement of the quantitative Subspace theorem},
Compositio Math. 101 (1996), 225--311.

\sm \item{[17]}
S. Ferenczi,
{\it Rank and symbolic complexity},
Ergodic Th. Dyn. Systems 16 (1996), 663--682.

\sm \item{[18]}
 S. Ferenczi \& C. Mauduit,
 {\it Transcendence of numbers with a low complexity expansion},
 J.~ Number~ Theory 67 (1997), 146--161.
 
\sm \item{[19]}
 J. Hartmanis \& R. E. Stearns,
{\it On the computational complexity of algorithms},
Trans. Amer. Math. Soc. 117 (1965), 285--306.

\sm \item{[20]}
 A. J. Kempner,
{\it On Transcendental Numbers},
Trans. Amer. Math. Soc. 17 (1916), 476--482.

\sm \item{[21]}
J. H. Loxton \& A. J. van der Poorten,
 {\it Arithmetic properties of the solutions of a class of
 functional equations},
 J. Reine Angew. Math. 330 (1982), 159--172.
 
\sm \item{[22]}
J. H. Loxton \& A. J. van der Poorten,
 {\it Arithmetic properties of automata: regular sequences},
 J. Reine Angew. Math. 392 (1988), 57--69.

\sm \item{[23]}
 K. Mahler,
 {\it Arithmetische Eigenschaften der L\"osungen einer Klasse von 
Funktionalgleichungen},
Math. Annalen, 101 (1929), 342--366. Corrigendum 103 (1930), 532.
  
\sm \item{[24]}
K. Mahler,
{\it Arithmetische Eigenschaften einer Klasse 
transzendental-transzendenter Funktionen}, Math. Z., 32 (1930), 545-585.

\sm \item{[25]}
 K. Mahler,
 {\it \"Uber das Verschwinden von Potenzreihen mehrerer
 Ver\"anderlichen in speziellen Punktfolgen},
Math. Ann., 103 (1930), 573--587.

\sm \item{[26]}
M. Morse \& G. A. Hedlund,
{\it Symbolic dynamics},
Amer. J. Math. 60 (1938), 815--866.

\sm \item{[27]}
M. Morse \& G. A. Hedlund,
{\it Symbolic dynamics II},
Amer. J. Math. 62 (1940), 1--42.
   
\sm \item{[28]}
P. Morton \& W. J. Mourant, 
{\it Digit patterns and transcendental numbers}, 
J. Austral. Math. Soc. Ser. A 51 (1991), 216--236.

\sm \item{[29]}
Ku. Nishioka,
Mahler functions and transcendence,
Lecture Notes in Math. 1631, Springer-Verlag, Berlin, 1996.

\sm \item{[30]}
A. Renyi,
{\it Representations for real numbers and their ergodic
properties},
Acta Math. Sci. Hungar. 8 (1957), 477--493.

\sm \item{[31]}
D. Ridout,
{\it Rational approximations to algebraic numbers},
Mathematika 4 (1957), 125--131.

\sm \item{[32]}
R. N. Risley \& L. Q. Zamboni,
{\it A generalization of Sturmian sequences: combinatorial
              structure and transcendence},
Acta Arith. 95 (2000), 167--184.
  
\sm \item{[33]}
H. P. Schlickewei, 
{\it On products of special linear forms with algebraic coefficients},
Acta Arith. 31 (1976), 389--398.

\sm \item{[34]}
H. P. Schlickewei, 
{\it The ${\goth p}$-adic Thue-Siegel-Roth-Schmidt theorem},
Arch. Math. (Ba\-sel) 29 (1977), 267--270.

\sm \item{[35]}
H. P. Schlickewei, 
{\it The quantitative Subspace Theorem for number fields},
Compositio Math. 82 (1992), 245--273.

\sm \item{[36]}
K. Schmidt, 
{\it On periodic expansions of Pisot numbers and Salem numbers},
Bull. London Math. Soc. 12 (1980), 269--278.

\sm \item{[37]}
 W. M. Schmidt,
{\it Norm form equations},
Ann. of Math. 96 (1972), 526--551.

\sm \item{[38]}
W. M. Schmidt, 
Diophantine approximation,
Lecture Notes in Mathematics 785 (1980) Springer.

\sm \item{[39]}
A. M. Turing,
{\it On computable numbers, with an application to the 
Entscheidungsproblem},
Proc. London Math. Soc. 42 (1937), 230--265.

\sm \item{[40]}
M. Waldschmidt,
{\it Un demi-si\`ecle de transcendance}. In:
Development of mathematics 1950--2000, pp. 1121--1186,
Birkh\"auser, Basel, 2000.

\vskip12mm

\noindent Boris Adamczewski   \hfill{Yann Bugeaud}

\noindent   CNRS, Institut Camille Jordan  
\hfill{Universit\'e Louis Pasteur}

\noindent   Universit\'e Claude Bernard Lyon 1 
\hfill{U. F. R. de math\'ematiques}

\noindent   B\^at. Braconnier, 21 avenue Claude Bernard
 \hfill{7, rue Ren\'e Descartes}

\noindent   69622 VILLEURBANNE Cedex (FRANCE)   
\hfill{67084 STRASBOURG Cedex (FRANCE)}

\vskip2mm
 
\noindent {\tt Boris.Adamczewski@math.univ-lyon1.fr}
\hfill{{\tt bugeaud@math.u-strasbg.fr}}

\bye